\documentclass[11pt]{amsart}
\usepackage{amsmath,amssymb,amsthm}
\usepackage{accents,graphicx,tikz}
\usepackage[shortlabels]{enumitem}
\newcommand{\longsim}{\scalebox{1.3}[0.8]{\ensuremath{\sim}}}
\newcommand{\MG}{\underaccent{\longsim}{\mathcal{M}}}
\newcommand{\set}[1]{\{#1\}}
\newcommand{\Z}{\mathbb{Z}}

\theoremstyle{definition}
\newtheorem{thm}{Theorem}[section]
\newtheorem{lem}[thm]{Lemma}

\newtheorem{defn}[thm]{Definition}
\newtheorem{model}[thm]{Model}
\theoremstyle{remark}
\newtheorem{rem}[thm]{Remark}

\title[Finite models for the order axioms]{Order Axioms of the Hilbert Plane:\\ Redundancies and Finite Independence Models}
\author{Roberto Volpe}
\date{\today}

\begin{document}
\begin{abstract}
Results on the redundancy and on the independence of the order axioms of the Hilbert plane are scattered in the literature, from Veblen's 1904 remark that the independence of the axioms of order can be shown by finite models to the redundancies removed from Hilbert's axioms of order in later editions of the \emph{Grundlagen}. The aim of this note is to collect them in an organic and self-contained form for the incidence and order axioms as presented by Greenberg. In doing so we single out three redundancies inside the order axioms: the symmetry of betweenness, the existence part of trichotomy, and the second clause of plane separation. Each of them dictates which part of the corresponding axiom an independence model must violate, and we give explicit finite models of three, three, nine and twenty-five points. A model of all the incidence and order axioms, on the contrary, is necessarily infinite, so the independence of any axiom added to them requires infinite models; for two of the order axioms the finite models also give independence from all the remaining axioms of the Hilbert plane.
\end{abstract}
\maketitle

\section{Introduction}

The order axioms of the Hilbert plane are usually introduced in textbooks with an informal remark on infinity. Greenberg asks the reader to show that, once the betweenness axioms are added to the incidence axioms, every line has at least five points, and to
\begin{quote}
``give an informal argument to show that every segment (a fortiori, every line) has an infinite number of points lying on it (a formal proof requires the technique of mathematical induction)'' \cite[Chapter~3, Exercise~8]{Gre}.
\end{quote}
Hartshorne states the same fact as an exercise \cite[Exercise~7.4]{Har}, and in the next exercise asks for a model of betweenness on the points of a single line with only finitely many points, suggesting the integers modulo~$5$ with $b$ between $a$ and $c$ when $b$ is their midpoint \cite[Exercise~7.5]{Har}.

The phenomenon behind these exercises was pointed out by Veblen in 1904 for his own axiom system:
\begin{quote}
``The independence proofs of the first eight axioms can all be exhibited by means of finite classes of elements or points. That such independence proofs are not possible for IX--XII is a consequence of theorem~12, based on I--VIII, that every line contains an infinitude of points'' \cite[\S7]{Veb}.
\end{quote}
Compared with our presentation, which follows Greenberg, Veblen's system takes order as the only primitive relation, defines lines, and uses Pasch's axiom in place of plane separation.

The aim of this note is to collect in an organic and self-contained form results on the incidence and order axioms, as presented by Greenberg, which are scattered in the literature. The main fact is that the independence of each order axiom from the incidence axioms and from the remaining order axioms can be proved with finite models only. A model of all the incidence and order axioms, on the contrary, is necessarily infinite \cite[Exercise~7.4]{Har}; as a consequence, the independence of any axiom added to them requires infinite models.

Three parts of the order axioms are redundant:
\begin{enumerate}[(a)]
\item the symmetry of betweenness follows from the incidence axioms, the second and fourth order axioms and the first part of the first order axiom;
\item the existence part of trichotomy follows from the incidence axioms, the first, second and fourth order axioms and the uniqueness part of trichotomy;
\item the second clause of plane separation follows from the incidence axioms, the first and third order axioms and the first clause of plane separation.
\end{enumerate}
Each redundancy dictates which part of the corresponding axiom an independence model must violate, and the finite models of Section~\ref{sec:models} violate exactly those parts.

The axioms of order in Hilbert's \emph{Grundlagen} were revised through several editions as redundancies were discovered \cite{SF,BSV}; the fact that a line cannot meet the three sides of a triangle in interior points is attributed to Pasch by Hallett and Majer in an editorial note to Hilbert's 1898--99 lectures \cite[Chapter~4, note~36]{HM};  Scott and Fleuriot report that by the tenth edition the existence part of trichotomy was derived as a theorem, with a proof attributed to Wald, and that a further order axiom was removed with a proof due to E.~H.~Moore \cite{SF}; the existence part of trichotomy is also derived from the other axioms in \cite[Theorem~8]{BSV}.
Our purpose is to present them together, for Greenberg's axioms, with complete proofs and explicit models, and to show how the redundancies and the finite models fit together.

The paper is organised as follows. Section~\ref{sec:axioms} recalls the axioms and the definitions of segment and side. Section~\ref{sec:red} proves the three redundancies. Section~\ref{sec:models} presents the four finite models. Section~\ref{sec:congr} remarks that two of them also give independence from the congruence axioms, and Section~\ref{sec:concl} concludes.

\section{Axioms and definitions}\label{sec:axioms}

A \emph{Hilbert plane} is a model of the incidence, order and congruence axioms \cite[\S10]{Har}. We use Greenberg's axioms \cite[Chapters~2--3]{Gre} with the notation of \cite{Vol1}, where the system is denoted $\MG$. The primitive terms are \emph{point}, \emph{line}, \emph{incidence}, the ternary relation of \emph{betweenness}, written $A-B-C$, and the relations of \emph{congruence} of segments and of angles, both written $\equiv$. Sections~\ref{sec:red} and~\ref{sec:models} use only the incidence and order axioms.

\begin{enumerate}[label=\textrm{[}\textbf{I}\arabic*\textrm{]}]
  \item\label{I1} For every point $P$ and for every point $Q\neq P$, there exists a unique line $l$ incident with $P$ and $Q$.
  \item\label{I2} For every line $l$ there exist at least two distinct points that are incident with $l$.
  \item\label{I3} There exist three distinct points with the property that no line is incident with all three of them.
\end{enumerate}
\begin{enumerate}[label=\textrm{[}\textbf{O}\arabic*\textrm{]}]
  \item\label{O1} If $A-B-C$, then $A$, $B$, and $C$ are three distinct collinear points, and $C-B-A$.
  \item\label{O2} Given any two distinct points $B$ and $D$, there exist points $A$, $C$, and $E$ lying on $\overleftrightarrow{BD}$ such that $A-B-D$, $B-C-D$, and $B-D-E$.
  \item\label{O3} If $A$, $B$, and $C$ are three distinct points lying on the same line, then one and only one of the points is between the other two.
  \item\label{O4} For every line $l$ and for any three points $A$, $B$, and $C$ not lying on $l$:
  \begin{enumerate}[label=(\roman*)]
    \item if $A$ and $B$ are on the same side of $l$ and $B$ and $C$ are on the same side of $l$, then $A$ and $C$ are on the same side of $l$;
    \item if $A$ and $B$ are on opposite sides of $l$ and $B$ and $C$ are on opposite sides of $l$, then $A$ and $C$ are on the same side of $l$.
  \end{enumerate}
\end{enumerate}
\begin{enumerate}[label=\textrm{[}\textbf{C}\arabic*\textrm{]}]
  \item\label{C1} If $A$ and $B$ are distinct points and if $A'$ is any point, then for each ray $r$ emanating from $A'$ there is a unique point $B'$ on $r$ such that $B'\neq A'$ and $AB\equiv A'B'$.
  \item\label{C2} If $AB\equiv CD$ and $AB\equiv EF$, then $CD\equiv EF$. Moreover, every segment is congruent to itself.
  \item\label{C3} If $A-B-C$, $A'-B'-C'$, $AB\equiv A'B'$, and $BC\equiv B'C'$, then $AC\equiv A'C'$.
  \item\label{C4} Given an angle $\angle BAC$ (where, by definition of ``angle'', $\overrightarrow{AB}$ is not opposite to $\overrightarrow{AC}$), and given any ray $\overrightarrow{A'B'}$ emanating from a point $A'$, then there is a unique ray $\overrightarrow{A'C'}$ on a given side of line $\overleftrightarrow{A'B'}$ such that $\angle B'A'C'\equiv\angle BAC$.
  \item\label{C5} If $\angle A\equiv\angle B$ and $\angle A\equiv\angle C$, then $\angle B\equiv\angle C$. Moreover, every angle is congruent to itself.
  \item\label{C6} If two sides and the included angle of one triangle are congruent respectively to two sides and the included angle of another triangle, then the two triangles are congruent.
\end{enumerate}
We call \emph{first part} of \ref{O1} the statement that $A$, $B$, $C$ are distinct and collinear, and \emph{symmetry} the statement $C-B-A$; \emph{existence} and \emph{uniqueness} in \ref{O3} are the two halves of ``one and only one''.

\begin{defn}[{\cite[Chapter~3]{Gre}}]
For distinct points $A$, $B$, the \emph{segment} $AB$ is the set formed by $A$, $B$ and the points $P$ with $A-P-B$. Two points $A$, $B$ not on a line $l$ are on the \emph{same side} of $l$ if $A=B$ or the segment $AB$ does not meet $l$; they are on \emph{opposite sides} of $l$ if $A\neq B$ and $AB$ meets $l$.
\end{defn}

Note that the axioms do not require the points in \ref{O2} to be distinct from each other, nor the three points in \ref{O4} to be distinct; the finite models of Section~\ref{sec:models} use this freedom.

We shall use two consequences of \ref{O4}, indicating the axioms used in their proofs: Lemma~\ref{lem:iii} is stated without proof in \cite[Chapter~3, p.~77]{Gre}, and we give its proof; Lemma~\ref{lem:pasch} is proved in \cite[Chapter~3]{Gre}.

\begin{lem}[{\cite[Chapter~3, p.~77]{Gre}}]\label{lem:iii}
Let $l$ be a line and $A$, $B$, $C$ points not on $l$. If $A$ and $B$ are on opposite sides of $l$ and $B$ and $C$ are on the same side of $l$, then $A$ and $C$ are on opposite sides of $l$.
\end{lem}
\begin{proof}
The proof uses only \ref{O4}(i) and the symmetry in \ref{O1}; in particular it does not use \ref{O4}(ii), so the lemma can be used in the proof of Theorem~\ref{thm:sep} without circularity.

If $A$ and $C$ were on the same side of $l$, then by \ref{O4}(i) applied to $A$, $C$, $B$ the points $A$ and $B$ would be on the same side of $l$; ``same side'' and ``opposite sides'' exclude each other by definition. The step from ``$B$ and $C$ on the same side'' to ``$C$ and $B$ on the same side'' uses that the segment $BC$ coincides with $CB$, which follows from the symmetry in \ref{O1}.
\end{proof}

\begin{lem}[Pasch, {\cite[Pasch's Theorem]{Gre}}]\label{lem:pasch}
If $A$, $B$, $C$ are noncollinear and a line $l$ not through $A$, $B$, $C$ meets the segment $AB$, then $l$ meets $AC$ or $BC$.
\end{lem}
\begin{proof}
For the proof see \cite[Pasch's Theorem]{Gre}. We only point out that it uses only [\textbf{I}], \ref{O1} and \ref{O4}; in particular it does not use \ref{O3}, so the lemma can be used in the proof of Theorem~\ref{thm:ex} without circularity.
\end{proof}

\section{Three redundancies}\label{sec:red}
This section proves the three redundancies (a)–(c) of the Introduction.

\begin{thm}\label{thm:sym}
The symmetry in \ref{O1} follows from [\textbf{I}], \ref{O2}, \ref{O4} and the first part of \ref{O1}.
\end{thm}
\begin{proof}
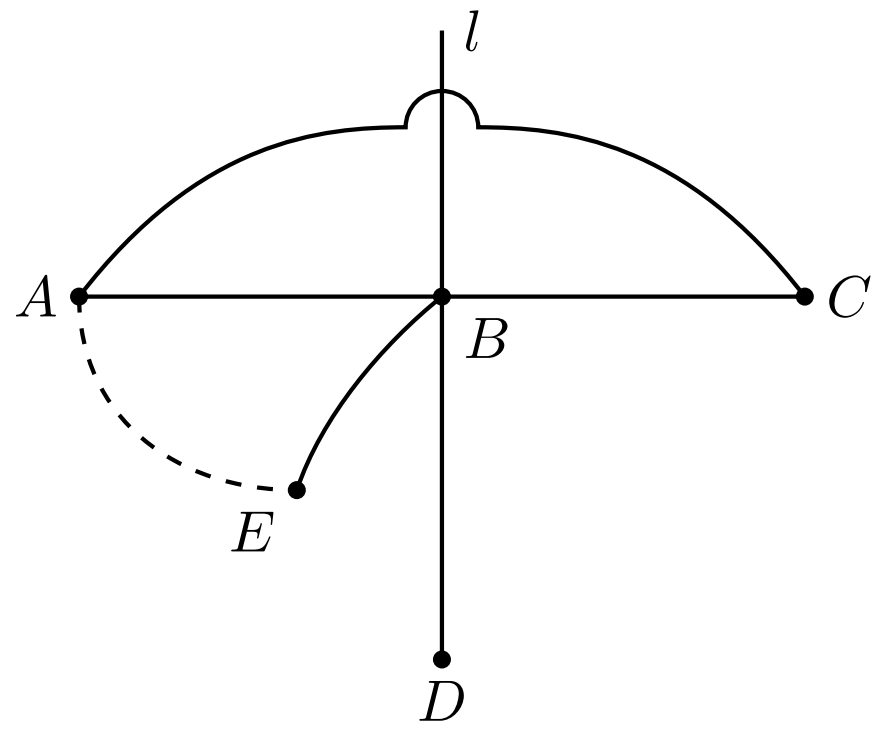
\begin{figure}
  \centering
\begin{tikzpicture}[scale=0.8,thick,
  pto/.style={circle,fill,inner sep=1.2pt}]
\coordinate (A) at (-3,0);
\coordinate (B) at (0,0);
\coordinate (C) at (3,0);
\coordinate (D) at (0,-3);
\coordinate (E) at (-1.2,-1.6);
\draw (0,2.2) -- (D);
\node[right] at (0,2.2) {$l$};
\draw (A) -- (C);
\draw (C) .. controls (2,1.3) and (1,1.4) .. (0.3,1.4)
      arc[start angle=0,end angle=180,radius=0.3]
      .. controls (-1,1.4) and (-2,1.3) .. (A);
\draw (B) .. controls (-0.5,-0.4) and (-1,-1) .. (E);
\draw[dashed] (A) .. controls (-3,-1) and (-2.2,-1.6) .. (E);
\node[pto] at (A) {}; \node[left] at (A) {$A$};
\node[pto] at (B) {}; \node[below right] at (B) {$B$};
\node[pto] at (C) {}; \node[right] at (C) {$C$};
\node[pto] at (D) {}; \node[below] at (D) {$D$};
\node[pto] at (E) {}; \node[below left] at (E) {$E$};
\end{tikzpicture}
\caption{Scheme of the proof of Theorem~\ref{thm:sym}: segments are drawn as paths; the path $CA+AE$ does not cross $l$, while the path $AC+CE$ crosses $l$ at $B$.}\label{fig:sym}
\end{figure}
Assume $A-B-C$ and, for contradiction, that $C-B-A$ does not hold.
\begin{enumerate}[label=(\arabic*)]
\item By \ref{I3} there is a point $D$ not on $\overleftrightarrow{AC}$; the line $l=\overleftrightarrow{BD}$ meets $\overleftrightarrow{AC}$ in $B$.
\item Since $B\in AC$, the points $A$ and $C$ are on opposite sides of $l$.
\item Since $C-B-A$ does not hold, $B\notin CA$. By the first part of \ref{O1} the segment $CA$ lies on $\overleftrightarrow{AC}$, which meets $l$ in $B$; hence $CA$ does not meet $l$, and $C$ and $A$ are on the same side of $l$.
\item By \ref{O2} applied to $C$ and $B$ there is $E$ with $C-B-E$; $E$ is distinct from $B$ and lies on $\overleftrightarrow{AC}$, hence not on $l$, and $B\in CE$ gives that $C$ and $E$ are on opposite sides of $l$.
\item By \ref{O4}(ii), from (2) and (4), $A$ and $E$ are on the same side of $l$.
\item By \ref{O4}(i), from (3) and (5), $C$ and $E$ are on the same side of $l$, contradicting (4).\qedhere
\end{enumerate}
\end{proof}

\begin{thm}\label{thm:ex}
The existence part of \ref{O3} follows from [\textbf{I}], \ref{O1}, \ref{O2}, \ref{O4} and the uniqueness part of \ref{O3}.
\end{thm}
\begin{proof}
\begin{figure}
\centering
\begin{tikzpicture}[scale=1.2,thick,
  pto/.style={circle,fill,inner sep=1.2pt}]
\coordinate (A) at (0,0);
\coordinate (B) at (2.5,0);
\coordinate (C) at (4.5,0);
\coordinate (D) at (2,1.6);
\coordinate (E) at (1.5,3.2);
\coordinate (F) at (2.571,2.057);
\coordinate (G) at (1.038,2.215);
\draw (-0.6,0) -- (5.1,0) node[right] {$m$};
\draw (E) -- (B);
\draw[dotted] (E) -- (C);
\draw[dotted] (A) -- (E);
\draw[dashed] (A) -- (F);
\draw[dashed] (C) -- (G);
\node[pto] at (A) {}; \node[below] at (A) {$A$};
\node[pto] at (B) {}; \node[below] at (B) {$B$};
\node[pto] at (C) {}; \node[below] at (C) {$C$};
\node[pto] at (D) {}; \node[right=3pt] at (D) {$D$};
\node[pto] at (E) {}; \node[above] at (E) {$E$};
\node[pto] at (F) {}; \node[right] at (F) {$F$};
\node[pto] at (G) {}; \node[left] at (G) {$G$};
\end{tikzpicture}
\caption{Construction used in the proof of Theorem~\ref{thm:ex}; the fact that $A-B-C$ appears in the figure is to be read as the conclusion reached by following the construction.}\label{fig:ex}
\end{figure}
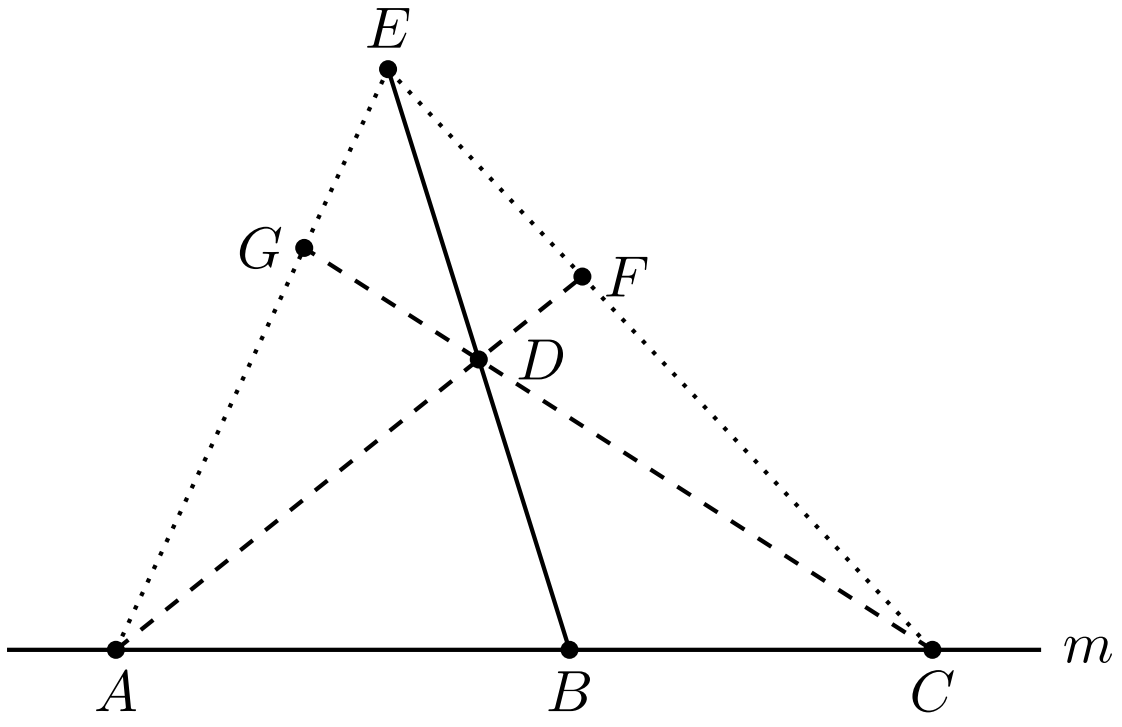
The proof is that of \cite[Theorem~8]{BSV}; the sub-items are ours and show where the hypotheses are used. By Lemma~\ref{lem:pasch}, Pasch's theorem uses only premises of the present theorem and does not use \ref{O3}, so it can be used without circularity.
\begin{enumerate}[label=(\arabic*)]
\item Let $A$, $B$, $C$ be distinct points of a line $m=\overleftrightarrow{AC}$, and assume that neither $A-C-B$ nor $B-A-C$ holds; we show that $A-B-C$ holds.
\item By \ref{I3} there is a point $D$ not on $\overleftrightarrow{AC}$.
\item By \ref{O2} there is a point $E$ on $\overleftrightarrow{DB}$ with $E-D-B$.
\item Applying Lemma~\ref{lem:pasch} to the triangle $BEC$ and the secant line $\overleftrightarrow{AD}$, we obtain a point $F$, the intersection of $\overleftrightarrow{AD}$ with $EC$, with $E-F-C$.
\begin{itemize}
  \item Since $B-A-C$ does not hold, $A$ is not in $BC$, and the line $\overleftrightarrow{AD}$ must meet $EC$ in the point $F$.
\end{itemize}
\item In the same way there is a point $G$, the intersection of $\overleftrightarrow{CD}$ with $AE$, with $A-G-E$.
\begin{itemize}
  \item Since $A-C-B$ does not hold, $C$ is not in $AB$, and the line $\overleftrightarrow{CD}$ must meet $AE$ in the point $G$.
\end{itemize}
\item The line $\overleftrightarrow{CG}$ is a secant line for the triangle $AEF$, hence $A-D-F$.
\begin{itemize}
  \item Since $E-F-C$, by the uniqueness in \ref{O3} the point $C$ is not in $EF$, and $\overleftrightarrow{CD}$ meets $AF$ in the point $D$.
\end{itemize}
\item Considering the triangle $AFC$ and the secant line $\overleftrightarrow{DE}$, it follows that $A-B-C$.
\begin{itemize}
  \item Since $E-F-C$, by the uniqueness in \ref{O3} the point $E$ is not in $FC$, and $\overleftrightarrow{DE}$ meets $AC$ in the point $B$.\qedhere
\end{itemize}
\end{enumerate}
\end{proof}

\begin{thm}\label{thm:sep}
Clause (ii) of \ref{O4} follows from [\textbf{I}], \ref{O1}, \ref{O3} and clause (i) of \ref{O4}.
\end{thm}
\begin{proof}
\begin{figure}
\centering
\begin{tikzpicture}[scale=0.70,thick,
  pto/.style={circle,fill,inner sep=1.2pt}]
\coordinate (A) at (-2,1.5);
\coordinate (B) at (0,-1.2);
\coordinate (C) at (2,1.6);
\coordinate (P) at (-0.889,0);
\coordinate (Q) at (0.857,0);
\coordinate (R) at (0.3,1.5575);
\draw (A) -- (B);
\draw (B) -- (C);
\draw[dashed] (A) -- (C);
\draw (-3,0) -- (3.2,0)
      .. controls (4.6,0) and (4.6,3.4) .. (2.6,3.4)
      .. controls (0.9,3.4) and (0.3,2.6) .. (R)
      .. controls (0.3,1.1) .. (0.15,0.7);
\node[left] at (-3,0) {$l$};
\node[pto] at (A) {}; \node[above left] at (A) {$A$};
\node[pto] at (B) {}; \node[below] at (B) {$B$};
\node[pto] at (C) {}; \node[above right] at (C) {$C$};
\node[pto] at (P) {}; \node[below left] at (P) {$P$};
\node[pto] at (Q) {}; \node[below right] at (Q) {$Q$};
\node[pto] at (R) {}; \node[above right] at (R) {$R$};
\node at (0.5,-2.2) {(a)};
\end{tikzpicture}
\qquad
\begin{tikzpicture}[scale=0.70,thick,
  pto/.style={circle,fill,inner sep=1.2pt}]
\coordinate (A) at (-1.4,0);
\coordinate (P) at (0,0);
\coordinate (B) at (1.6,0);
\coordinate (C) at (-0.9,1.3);
\coordinate (E) at (1.0,2.0);
\draw (-3,0) node[left] {$m$} -- (A) -- (P) -- (B)
      .. controls (2.8,0) and (2.9,-2.0) .. (1.4,-1.75)
      .. controls (1.0,-1.6) and (0.8,-1.15) .. (0.55,-0.8)
      -- (P) -- (C) -- (-1.35,1.95);
\draw (-0.9,-1.8) node[below] {$l$} -- (1.2,2.4);
\draw[dash dot] (1.75,-0.5) -- (0.85,2.5) node[above] {$\overleftrightarrow{EB}$};
\node[pto] at (A) {}; \node[below] at (A) {$A$};
\node[pto] at (P) {}; \node[below left] at (P) {$P$};
\node[pto] at (B) {}; \node[above right] at (B) {$B$};
\node[pto] at (C) {}; \node[right] at (C) {$C$};
\node[pto] at (E) {}; \node[right] at (E) {$E$};
\node at (0,-3) {(b)};
\end{tikzpicture}
\caption{Scheme of the proof of Theorem~\ref{thm:sep}: the two cases for $A$, $B$, $C$, (a) noncollinear and (b) collinear.}\label{fig:sep}
\end{figure}
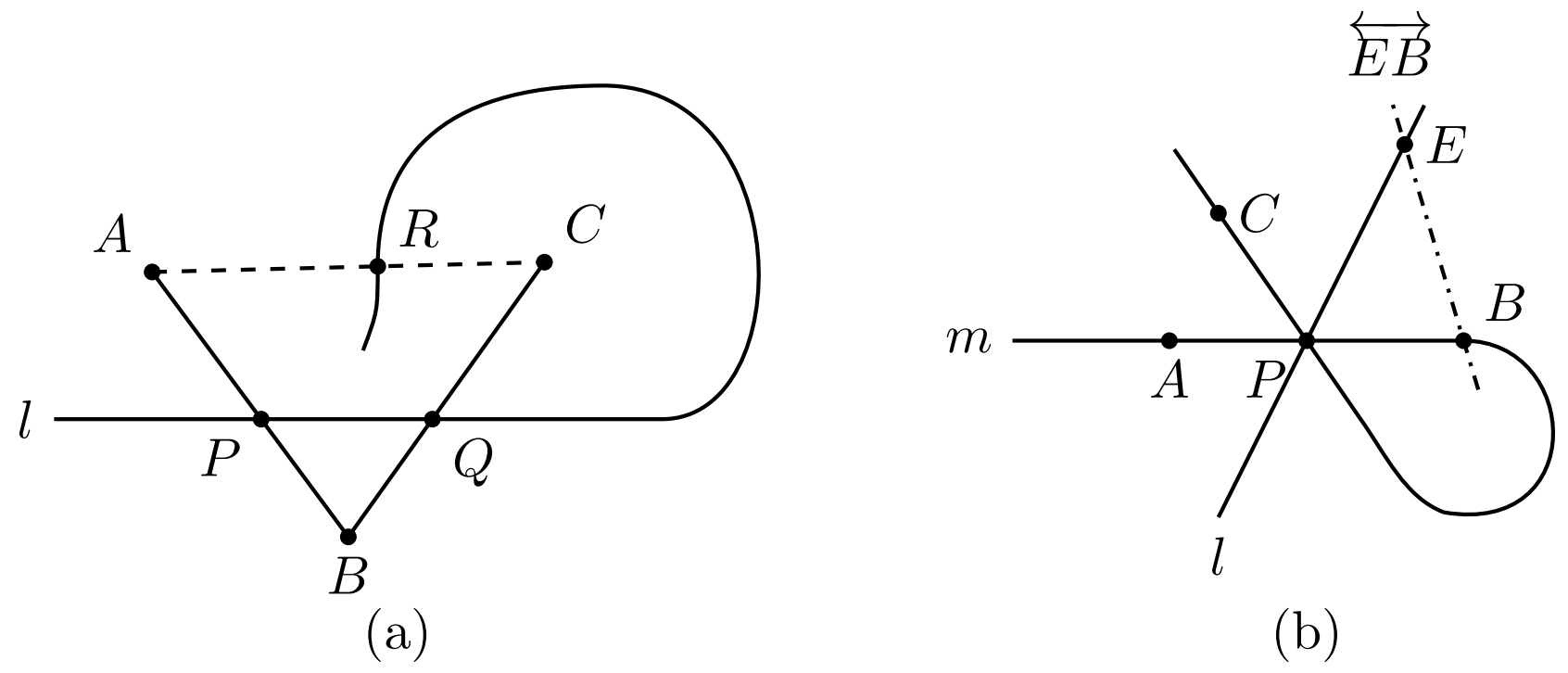
Let $l$ be a line and $A$, $B$, $C$ points not on $l$, with $A$, $B$ on opposite sides and $B$, $C$ on opposite sides of $l$. If $A=C$ there is nothing to prove; assume $A\neq C$ and, for contradiction, that $A$ and $C$ are on opposite sides of $l$. Then there are points $P$, $Q$, $R$ of $l$ with $A-P-B$, $B-Q-C$ and $A-R-C$.
\begin{enumerate}[label=(\arabic*)]
\item If $A$, $B$, $C$ are not collinear, then $P$, $Q$, $R$ are distinct: if for instance $P=Q$, the lines $\overleftrightarrow{AB}$ and $\overleftrightarrow{BC}$ would share the distinct points $P$ and $B$, and $A$, $B$, $C$ would be collinear. Since $P$, $Q$, $R$ lie on $l$, by \ref{O3} one of them is between the other two; since the three pairs of vertices are on opposite sides of $l$, after renaming the vertices we may assume $P-Q-R$ (Figure~\ref{fig:sep}(a)).
\item The points $P$ and $R$ are on opposite sides of $\overleftrightarrow{BC}$. Indeed neither lies on $\overleftrightarrow{BC}$, otherwise $A$, $B$, $C$ would be collinear; and, by $P-Q-R$, the segment $PR$ meets $\overleftrightarrow{BC}$.
\item The points $A$ and $P$ are on the same side of $\overleftrightarrow{BC}$. Indeed the segment $AB$ and the line $\overleftrightarrow{BC}$ meet only in $B$; if $A$ and $P$ were on opposite sides, $AP$ would meet $\overleftrightarrow{BC}$ in $B$, giving $A-B-P$, which contradicts the uniqueness in \ref{O3} since $A-P-B$.
\item By Lemma~\ref{lem:iii} applied to (2) and (3), $A$ and $R$ are on opposite sides of $\overleftrightarrow{BC}$. Since the segment $AC$ and the line $\overleftrightarrow{BC}$ meet only in $C$, $AR$ would meet $\overleftrightarrow{BC}$ in $C$, giving $A-C-R$, which together with $A-R-C$ contradicts \ref{O3}. Hence $A$, $B$, $C$ cannot be noncollinear.
\item Let $A$, $B$, $C$ lie on a line $m$. Since $A$, $B$, $C$ are not on $l$, $m\neq l$; hence $m$ and $l$ meet in a single point, so $P=Q=R$, and $A-P-B$, $B-P-C$, $A-P-C$. By \ref{O3} one of $A$, $B$, $C$ is between the other two; after renaming we may assume $A-B-C$.
\item By \ref{I3} there is a point $E$ not on $m$; the line $\overleftrightarrow{EB}$ meets $m=\overleftrightarrow{AC}$ in $B$. Since $A-B-C$, the points $A$ and $C$ are on opposite sides of $\overleftrightarrow{EB}$.
\item The points $A$ and $P$ are on the same side of $\overleftrightarrow{EB}$: otherwise $AP$ would meet $\overleftrightarrow{EB}$ in $B$, giving $A-B-P$, which contradicts the uniqueness in \ref{O3} since $A-P-B$.
\item By Lemma~\ref{lem:iii} applied to (6) and (7), $C$ and $P$ are on opposite sides of $\overleftrightarrow{EB}$, so $PC$ would meet $\overleftrightarrow{EB}$ in $B$, giving $P-B-C$, which together with $B-P-C$ contradicts \ref{O3}. Hence $A$, $B$, $C$ cannot be collinear.
\end{enumerate}
In both cases we reach a contradiction, so $A$ and $C$ are on the same side of $l$.
\end{proof}

\begin{rem}\label{rem:constraints}
By Theorems~\ref{thm:sym}--\ref{thm:sep}, a model showing the independence of \ref{O1} from the other incidence and order axioms must violate the first part of \ref{O1}; a model for \ref{O3} must violate uniqueness; a model for \ref{O4} must violate clause (i).
\end{rem}

\section{Finite independence models}\label{sec:models}

\begin{model}[\ref{O2} fails]\label{mod:O2}
Points $A$, $B$, $C$; lines $\set{A,B}$, $\set{B,C}$, $\set{C,A}$; betweenness empty. Then [\textbf{I}] hold; \ref{O1} and \ref{O3} hold vacuously, since no line has three points; every segment consists of its endpoints, outside each line there is exactly one point, so \ref{O4}(i) holds (a point is on the same side as itself) and \ref{O4}(ii) holds vacuously. \ref{O2} fails, since a line contains only two points.
\end{model}

\begin{model}[\ref{O1} fails]\label{mod:O1}
Points and lines as in Model~\ref{mod:O2}; betweenness
\begin{equation*}
X-Y-Z\quad\Longleftrightarrow\quad X=Y\ \text{or}\ Y=Z.
\end{equation*}
For $X\neq Y$ the segment $XY$ is still $\set{X,Y}$, so segments and sides are those of Model~\ref{mod:O2}, and [\textbf{I}], \ref{O3}, \ref{O4} hold as there. \ref{O2} holds: given $B\neq D$, the points $B$, $B$, $D$ lie on $\overleftrightarrow{BD}$ and $B-B-D$, $B-B-D$, $B-D-D$. \ref{O1} fails, since $A-A-B$ holds with $A$, $A$, $B$ not distinct. As required by Remark~\ref{rem:constraints}, it is the first part of \ref{O1} that fails.
\end{model}

\begin{model}[\ref{O3} fails]\label{mod:O3}
The affine plane of order~$3$: points $(x,y)\in\Z_{3}^{2}$, lines the solution sets of $ax+by=c$ modulo~$3$ with $(a,b)\neq(0,0)$; each line has three points. Betweenness holds for every triple of distinct points of a line, in every order. Then [\textbf{I}] hold, \ref{O1} holds by construction, \ref{O2} holds by taking the three required points equal to the third point of $\overleftrightarrow{BD}$. Every segment is a whole line, so two points not on $l$ are on the same side of $l$ exactly when they lie on the same parallel to $l$; there are two parallels to $l$ besides $l$, so both clauses of \ref{O4} hold. \ref{O3} fails, because each point of a line is between the other two: uniqueness fails while existence holds.
\end{model}

\begin{model}[\ref{O4} fails]\label{mod:O4}
The affine plane of order~$5$: points $(x,y)\in\Z_{5}^{2}$, lines the solution sets of $ax+by=c$ modulo~$5$ with $(a,b)\neq(0,0)$; each line has five points. Betweenness is the midpoint relation:
\begin{equation*}
A-B-C\quad\Longleftrightarrow\quad A,\ B,\ C\ \text{distinct and}\ (A+C)/2=B,
\end{equation*}
the operations being taken on coordinates modulo~$5$; on a single line this is Hartshorne's example \cite[Exercise~7.5]{Har}, and it is also the order on the lines of five points of Veblen's twenty-one-point model for his axiom of Pasch \cite[\S\S7--8]{Veb}. For $P\neq Q$ the segment $PQ$ consists of $P$, $Q$ and their midpoint $M$. If $l$ has equation $f(x,y)=0$ with $f(x,y)=ax+by-c$, then $f(M)=(f(P)+f(Q))/2$; hence two points $P$, $Q$ not on $l$ are on opposite sides of $l$ if and only if $P\neq Q$ and $f(P)+f(Q)=0$.
\begin{itemize}
\item [\textbf{I}] hold.
\item \ref{O1}: if $A-B-C$, then $B$ is the midpoint of $A$ and $C$, hence lies on $\overleftrightarrow{AC}$, and the condition is symmetric in $A$ and $C$.
\item \ref{O2}: for $B\neq D$, the points $2B-D$, $(B+D)/2$, $2D-B$ satisfy the three relations and are distinct from $B$ and $D$.
\item \ref{O3}: trichotomy says that exactly one of three distinct collinear points is the midpoint of the other two. If none were, the three midpoints would give six distinct points on a line of five points; if two were, two of the points would coincide.
\item \ref{O4}(i) fails: for $l\colon y=0$ and $A=(0,1)$, $B=(0,2)$, $C=(0,4)$ we have $f(A)+f(B)=3$, $f(B)+f(C)=1$, $f(A)+f(C)=0$.
\item \ref{O4}(ii) holds: for any line $r\colon g=0$ and points $X$, $Y$, $Z$ not on $r$, from $g(X)=-g(Y)$ and $g(Y)=-g(Z)$ we get $g(X)+g(Z)=2g(X)\neq0$.
\end{itemize}
As required by Remark~\ref{rem:constraints}, it is clause (i) of \ref{O4} that fails.
\end{model}

\section{A remark on the congruence axioms}\label{sec:congr}

Models~\ref{mod:O2} and~\ref{mod:O1} can be completed with congruence relations satisfying all the congruence axioms \ref{C1}--\ref{C6}: take as segment congruence and angle congruence the total relations on the three segments and on the three angles. In both models every ray contains exactly one point besides its origin and every line has exactly one point outside it, so the transport axioms \ref{C1} and \ref{C4} hold; \ref{C2} and \ref{C5} hold because total relations are equivalence relations; \ref{C3} and \ref{C6} hold because their conclusions are always true (in Model~\ref{mod:O1} the hypotheses of \ref{C3} are never met, since one of the two segments would have coinciding endpoints). Hence \ref{O1} and \ref{O2} are independent of all the remaining twelve axioms of $\MG$, by finite models.

\section{Conclusions}\label{sec:concl}

The independence of each order axiom from the incidence axioms and from the remaining order axioms can be shown with finite models only, and each redundancy of Section~\ref{sec:red} tells which part of an axiom such a model must violate: the first part of \ref{O1}, the uniqueness in \ref{O3}, clause (i) of \ref{O4}. A model of all the incidence and order axioms is necessarily infinite \cite[Exercise~7.4]{Har}, so infinite models enter only with the axioms added to them.

\end{document}